\pdfoutput=1
\documentclass{article}
\usepackage{amsthm,amsmath,amssymb,xspace,url,graphicx,caption,subcaption,a4wide}
\RequirePackage{natbib}
\RequirePackage[colorlinks,citecolor=blue,urlcolor=blue]{hyperref}

\newcommand{\R}{\mathbb{R}}

\newcommand{\N}{\mathbb{Z}_{\geq 0}}

\newcommand{\graf}{{\cal G}}
\newcommand{\E}{{\mbox{E}}}
\newcommand{\cd}{\,|\,}
\newcommand{\etal}{\mbox{\textit{et al.}}}

\newcommand{\bp}{P}

\newcommand{\bo}{\mbox{\boldmath $0$}}
\newcommand{\bx}{\mbox{\boldmath $x$}}
\newcommand{\bX}{\mbox{\boldmath $X$}}
 \newcommand{\cX}{{\cal
    X}}  

\newcommand{\iid}{{independent and identically distributed}\xspace}

\title{A Note on Bayesian Model Selection for Discrete Data Using
  Proper Scoring Rules}

\author{A.~Philip~Dawid\thanks{Department of Pure Mathematics and
    Mathematical Statistics, University of Cambridge, U.K.} \and
  {Monica~Musio}\thanks{Dipartiment of Mathematics, University of
    Cagliari, Italy} \and {Silvia~Columbu}\thanks{Dipartiment of
  Mathematics, University of Cagliari, Italy}}

\begin{document}
\maketitle
\begin{abstract}
  \noindent We consider the problem of choosing between parametric
  models for a discrete observable, taking a Bayesian approach in
  which the within-model prior distributions are allowed to be
  improper.  In order to avoid the ambiguity in the marginal
  likelihood function in such a case, we apply a homogeneous scoring
  rule.  For the particular case of distinguishing between Poisson and
  Negative Binomial models, we conduct simulations that indicate that,
  applied prequentially, the method will consistently select the true
  model.\\

\end{abstract}

\section{Introduction}
\label{sec:intro}
It is well known that Bayesian model selection with improper
within-model prior distributions is not well-defined, owing to the
presence of an arbitrary multiplicative constant in each term of the
marginal likelihood function.  Recently \citep{apd/mm:modsel} it has
been shown how this problem can be overcome if one replaces negative
log-likelihood (the {\em log score\/}) by another, homogeneous, proper
scoring rule \citep{mfp/apd/sll:plsr}---since then the arbitrary
constants do not enter into the formulae.  That paper considered the
case of continuous variables and, in particular, the Hyv\"arinen
scoring rule \citep{Hyvarinen:2005es}, and showed that this approach
will generally lead to consistent selection of the correct model.

The above approach can not be applied directly when the data are
discrete, since then we need to use scoring rules specifically adapted
to the discrete case, as characterised in \citet{PLSRdiscrete}.  Here we 
investigate, by example, such a discrete data problem.  In particular
we consider the problem of distinguishing between the Poisson and the
Negative Binomial distributions.  Simulations indicate that the method
will again deliver consistent selection of the true model.

\section{Local scoring rules}
\label{sec:loc}

Let $\cX$ be a discrete sample space endowed with a structure whereby
with each $x\in\cX$ is associated a {\em neighbourhood\/}
$N_x\subseteq\cX$, containing $x$.  In \citet{PLSRdiscrete} it was
shown how to define a {\em proper local scoring rule\/} $S(x,P)$ on
$\cX$, where $x\in{\cal X}$, and $P$ is a distribution over
${\cal X}$.  The rule is {\em proper\/} if, for all $P$,
$S(P,Q) := \E_{X\sim P}S(X,Q)$ is minimised for $Q = P$, and {\em
  local\/} if $S(x,P)$ depend on $P$ only through the probabilities it
assigns to points in $N_x$.  Under a condition on the neighbourhoods,
we can define an undirected graph $\graf$ on $\cX$ such that we can
take $y\in N_x$ just when $x$ and $y$ are identical or are adjacent in
$\graf$.  Then all proper local scorings can be characterised, and (on
excluding the log score, yielding what are termed {\em key local\/}
proper scoring rules) any of these will be {\em homogeneous\/} in the
sense that its value is unchanged when all probabilities in $N_x$ are
scaled by the same positive constant.

In particular, suppose the sample space ${\cal X}$ is the set $\N$ of
non-negative integers, and we regard $x$ and $y$ as neighbours if and
only if they differ by at most 1.  It is shown in \citet{PLSRdiscrete}
that a key local scoring rule adapted to this structure has the form
\begin{equation}
  \label{eq:pairrule}
  S(x, \bp) = G_{x-1}'\left\{\frac{p(x)}{p(x-1)}\right\} +
  G_x\left\{\frac{p(x+1)}{p(x)}\right\} -
  \frac{p(x+1)}{p(x)}\, G_x'\left\{\frac{p(x+1)}{p(x)}\right\}\quad(x=0,1,\ldots)
\end{equation}
where, for each $x\in\N$, $p(x) = P(X=x)$, $G_x$ is a concave function
on $\R^+$, and the first term in \eqref{eq:pairrule} is absent if
$x=0$.  It is clear from the way in which ratios enter
\eqref{eq:pairrule} that such a scoring rule is homogeneous.

The cumulative score \eqref{eq:pairrule} based on an \iid sample
$(x_1, \ldots, x_n)$ in which the frequency of $y$ is $f_y$
$(y=0,1,\ldots)$ is
\begin{equation}
  \label{eq:totalscore}
  \sum_{y=0}^\infty f_y G_y(v_y) + \left(f_{y+1}-f_y v_y\right) G_y'(v_y)
\end{equation}
with $v_y := p(y+1)/p(y)$.  If for example we wished to fit the
Poisson model $p(x) \propto \theta^x/x!$, we might estimate $\theta$
by minimising the total empirical score
\begin{equation}
  \label{eq:poissscore}
  \sum_{y=0}^\infty f_y G_y\left(\frac{\theta}{y+1}\right) +
  \left(f_{y+1}- \frac{f_y}{y+1}\,\theta\right) G_y'\left(\frac{\theta}{y+1}\right). 
\end{equation}

In the sequel we shall use the special case of \eqref{eq:pairrule}
with
\begin{equation}
  \label{eq:gform0}
  G_x(v) = - (x+1)^a v^{m}/m(m-1) \quad(m > 0, m\neq 1).
\end{equation}
This gives the scoring rule
\begin{equation}
  \label{eq:sform}
  S(x,P) =  \left\{
    \begin{array}[c]{lr}
      m^{-1}\left\{{p(1)}/{p(0)}\right\}^m & (x=0)\\
                                           &\quad\\
      \{m(m-1)\}^{-1}\left[(m-1)(x+1)^a\left\{{p(x+1)}/{p(x)}\right\}^m\right.\\
      {}\quad\quad\quad\quad\quad\quad\quad\left. - mx^a\left\{{p(x)}/{p(x-1)}\}\right)^{m-1}\right] & (x > 0).
    \end{array}
  \right.
\end{equation}

\section{Bayesian Model Selection}
\label{sec:modsel}

Let ${\cal M}$ be a finite or countable class of statistical models
for the same observable $X \in{\cal X}$.  Each $M\in{\cal M}$ is a
parametric family, with parameter $\theta_M \in \Theta_M$, a
$d_M$-dimensional Euclidean space; when $M$ obtains, with parameter
value $\theta_M$, then $X$ has distribution $P_{\theta_M}$, with
density function (probability mass function) $p_M( x \cd \theta_M)$.
Having observed data $X = x$, we wish to make inference about which
model $M\in{\cal M}$ (and possibly which parameter-value $\theta_M$)
actually generated the data.

The Bayesian approach assigns, within each model $M$, a prior
distribution $\Pi_M$, with density $\pi_M(\cdot)$ say, for its
parameter $\theta_M$.  The associated {\em predictive distribution\/}
$P_M$ of $X$ (given only the validity of model $M$, but no information
on its parameter) has density function
\begin{equation}
  \label{eq:preddens}
  p_M(x) = \int_{\Theta_M} p_M(x \cd \theta_M)\, \pi_M(\theta_M)\,d\theta_M.
\end{equation}
Any function over ${\cal M}$ proportional to $p_M(x)$ (considered as a
function of $M$, for fixed $x$) supplies the {\em marginal
  likelihood\/} function, $L(M)$, based on data $X=x$.  In typical
asymptotic scenarios, selection of the model maximising $L(M)$, or,
equivalently, minimising the {\em log score\/}
$S_L(x, P_M) := - \log p_M(x)$, will consistently select the true
model \citep{apd:postprob}.

``Objective'' Bayesian inference attempts to use standardised
within-model priors $\Pi_M$ intended to represent ``prior ignorance''.
In many applications, such an ``ignorance prior'' for $\theta_M$ is
not a genuine distribution, but rather an ``improper'' $\sigma$-finite
but not finite measure, with a ``density'' $\pi_M(\cdot)$ that does
not have a finite integral and so can not be normalised to be a proper
probability density.  Typically one writes
$\pi_M(\theta_M) \propto f_M(\theta_M)$, where $f_M$ is a given
non-integrable function and the constant of proportionality is not
specified.  Even without that specification, this allows mechanical
computation of a formal within-model-$M$ posterior density
$\pi_M(\theta_M \cd x)$, by application of Bayes's formula:
$\pi_M(\theta_M \cd x) \propto p_M(x \cd \theta_M)\,
\pi_M(\theta_M)\propto p_M(x \cd \theta_M)\, f_M(\theta_M)$.  This
will often yield an integrable function and hence the possibility of
normalisation to supply a genuine probability density.

However things do not work out so well when we turn to model
selection.  We have, for each model $M$,
\begin{displaymath}
\pi_M(\theta_M) = c_M f_M(\theta_M),   
\end{displaymath}
where $c_M$ is the unspecified proportionality constant.  This formally
leads to the marginal likelihood function
\begin{displaymath}
  L_M \propto c_M \, \int_{\Theta_M} p_M(x \cd \theta_M)\, f_M(\theta_M)\,d\theta_M.
\end{displaymath}
But since this involves the unspecified constants $c_M$, which could
vary arbitrarily with $M$, it is no longer meaningful to compare
models by means of their marginal likelihoods.

A way round this problem was proposed in \citet{apd/mm:modsel}: instead
of attempting to minimise the log score
$S_L(x, P_M) := - \log p_M(x)$, we replace that with another proper
scoring rule $S(x,P_M)$.  And if that scoring rule is homogeneous, it
will simply not involve the unspecified constant $c_M$.  In
\citet{apd/mm:modsel} a detailed analysis of this approach was
conducted for the case of continuous data and the Hyv\"arinen scoring
rule, and it was shown that it will typically deliver consistent
selection of the true model.

\section{Discrete model selection}
\label{sec:discsel}

We shall investigate empirically, for a simple example, the validity
of the above results when generalised to the case of discrete data.
We shall use the scoring rule \eqref{eq:sform}, and apply this to the
choice between a Poisson and a Negative Binomial model.  For this
purpose we first need to compute, for each of these models separately,
the appropriate score.

\section{Poisson model}
\label{sec:bpois}

Consider the Poisson model $X \sim {\cal P}(k\Lambda)$:
\begin{equation}
  \label{eq:pois}
  p(x \cd \lambda) = e^{-k\lambda}(k\lambda)^x/x! \quad(x=0,1,\ldots),
\end{equation}
with conjugate prior $\Lambda \sim \Gamma(\alpha,\beta)$:
\begin{equation}
  \label{eq:gamma}
  \pi(\lambda) = \frac{\beta^\alpha}{\Gamma(\alpha)}\lambda^{\alpha-1}e^{-\beta\lambda}.
\end{equation}
For propriety we require $\alpha>0$, $\beta>0$.

The predictive distribution $P$ has density function
\begin{equation}
  \label{eq:poispred}
  p(x) = \frac{\Gamma(\alpha+x)}{\Gamma(\alpha)x!}(1-\phi)^\alpha\phi^x
\end{equation}
with $\phi:= k/(\beta+k)$.

Then $p(x+1)/p(x) = \phi(x+\alpha)/(x+1)$, and so
\begin{eqnarray}
  \label{eq:poisform0}
  S(0,P) &=&  m^{-1}\alpha^m\phi^m\\
  \nonumber
  S(x,P) &= &\{m(m-1)\}^{-1}\left\{(m-1)\phi^m(x+1)^{a-m}(x+\alpha)^m\right.\\
  \label{eq:poisform} 
  &&{} - \left. m\phi^{m-1}x^{a-m+1}(x+\alpha-1)^{m-1}\right\}\quad(x>0).
\end{eqnarray}

 
\subsection{Multiple observations}
\label{sec:multiple}

Suppose now we have $N$ {\iid} observations $\bX_N = (X_1, \ldots,
X_N)$ from the above Poisson distribution.  We can apply the above
score in two different ways:
\begin{enumerate}
\item Apply direct to the sufficient statistic.
\item Apply prequentially to all observations.
\end{enumerate}

\subsubsection{Sufficient statistic}
\label{sec:suffpois}
The sufficient statistic is $T_N = \sum_{i=1}^NX_i$, with distribution
${\cal P}(Nk\Lambda)$.  So the score computed this way is simply
obtained from \eqref{eq:poisform0} and \eqref{eq:poisform} on replacing $x$
by $t_N$ and $k$ by $Nk$.  This gives
\begin{eqnarray}
  \label{eq:suffpoisform0}
  S_N(\bo,P) &=&  m^{-1}\alpha^m\phi_N^m\\
  \nonumber
  S_N(\bx,P) &= &\{m(m-1)\}^{-1}\left\{(m-1)\phi_N^m(t_N+1)^{a-m}(t_N+\alpha)^m\right.\\
  \label{eq:suffpoisform}
  &&{} - \left. m\phi_N^{m-1}t_N^{a-m+1}(t_N+\alpha-1)^{a-m+1}\right\}\quad(t_N>0)
\end{eqnarray}
where $\phi_N := Nk/(\beta+Nk)$.

\subsubsection{Prequential}
\label{sec:preqpois}

Now suppose we have already observed $\bX^{n-1} = \bx^{n-1}$.  The
posterior distribution of $\Lambda$ is 
\begin{displaymath}
  \Lambda \cd \bX^{n-1} = \bx^{n-1} \sim \Gamma\left\{\alpha+t_{n-1},\beta+(n-1)k\right\}. 
\end{displaymath}
So the predictive distribution of $X_n$, given the previous
observations $\bX^{n-1} = \bx^{n-1}$, is obtained from
\eqref{eq:poisform0} and \eqref{eq:poisform} on replacing $x$ with $x_n$,
$\alpha$ with $\alpha+t_{n-1}$, and $\beta$ with $\beta+(n-1)k$.  The
incremental contribution to the prequential score is thus given by:
\begin{eqnarray}
  \label{eq:preqpoisform0}
  S_n^*(0,P) &=&  m^{-1}(\phi^*_n)^m(\alpha+t_{n-1})^m\\
  \nonumber
  S_n^*(x_n,P) &= &\{m(m-1)\}^{-1}\left\{(m-1)(\phi_n^*)^m(x_n+1)^{a-m}(t_n+\alpha)^m\right.\\
  \label{eq:preqpoisform}
  &&{} - \left. m(\phi_n^*)^{m-1}x_n^{a-m+1}(t_n+\alpha-1)^{a-m+1}\right\}\quad(x_n>0)
\end{eqnarray}
with $\phi_n^* := k/(\beta+nk)$.

The total prequential score is obtained by summing this from $n=1$ to
$N$.

\subsection{Improper prior}

The usual improper prior is the formal limit with $\alpha, \beta
\downarrow 0$.  In this case \eqref{eq:suffpoisform0} and \eqref{eq:suffpoisform} become:
\begin{eqnarray}
  \label{eq:imp1suffpoisform0}
  S_N(\bo,P) &=& 0\\
  \nonumber
  S_N(\bx,P) &= &\{m(m-1)\}^{-1}\left\{(m-1)(t_N+1)^{a-m}t_N^m\right.\\
  \label{eq:imp1suffpoisform}
  &&{} - \left. mt_N^{a-m+1}(t_N-1)^{a-m+1}\right\}\quad(t_N>0).
\end{eqnarray}
Note that the score is well-defined even when all observations are
$0$, in which case the posterior is improper.

For the prequential version, we obtain, from \eqref{eq:preqpoisform0} and
\eqref{eq:preqpoisform}:
\begin{eqnarray}
  \label{eq:imp1preqpoisform0}
  S_n^*(0,P) &=&  t_n^m/mn^m  \\
  \nonumber
  S_n^*(x_n,P) &= &(x_n+1)^{a-m}t_n^m/mn^m \\
  \label{eq:imp1preqpoisform}
             &&{} - x_n^{a-m+1}(t_n-1)^{a-m+1}/(m-1)n^{m-1}\quad(x_n>0).
\end{eqnarray}
An alternative improper prior is the Jeffreys prior, having
$\alpha=1/2$, $\beta\downarrow 0$, which is easily handled similarly.

\section{Negative Binomial model}
\label{sec:negbin}

Now we consider an alternative model, the Negative Binomial $X \sim
{\cal NB}(s;\Theta)$, having
\begin{equation}
  \label{eq:negbin}
  p(x \cd \theta) =  \frac{(s+x-1)!}{x!(s-1)!}(1-\theta)^s\theta^x\quad(x=0,1,\ldots),
\end{equation}
with conjugate prior $\Theta \sim \beta(p,q)$:
\begin{equation}
  \label{eq:beta}
  \pi(\theta) = \frac{\Gamma(p+q)}{\Gamma(p)\Gamma(q)} \theta^{p-1}(1-\theta)^{q-1}.
\end{equation}
For propriety we require $p>0$, $q>0$.

The predictive density is
\begin{equation}
  \label{eq:negbinpred}
  p(x) = \frac{\Gamma(p+q)}{\Gamma(p)\Gamma(q)} \frac{(s+x-1)!}{x!(s-1)!}
  \frac{\Gamma(p+x)\Gamma(q+s)}{\Gamma(p+q+s+x)}.
\end{equation}
Then 
\begin{displaymath}
  \frac{p(x+1)}{p(x)} = \frac{(x+s)(x+p)}{(x+1)(x+p+q+s)},
\end{displaymath}
and so we have:
\begin{eqnarray}
  \label{eq:negbinomform0}
  S(0,P) &=&  m^{-1}(sp)^m(p+q+s)^{-m}\\
  \nonumber
  S(x,P) &=&\{m(m-1)\}^{-1}\left[(m-1)(x+1)^{a-m}\{(x+s)(x+p)\}^{m}(x+p+q+s)^{-m} \right.\\
  \label{eq:negbinomform}
  &&{} - \left. mx^{a-m+1}\{(x+s-1)(x+p-1)\}^{m-1}(x+p+q+s-1)^{-m+1}\right].
\end{eqnarray}


\subsection{Multiple observations}
\label{sec:multnegbin}
Again, we can handle multiple observations either by restricting to
the sufficient statistic, or by cumulating the prequential score.

\subsubsection{Sufficient statistic}
\label{sec:suffnegbinom}
The sufficient statistic is $T_N = \sum_{i=1}^NX_i$, with distribution
${\cal NB}(Ns, \Theta)$.  So the score computed this way is simply
obtained from \eqref{eq:negbinomform0} and \eqref{eq:negbinomform} on replacing $x$
by $t_N$ and $s$ by $Ns$.  This gives
\begin{eqnarray}
  \label{eq:suffnegbinomform0}
  S_N(\bo,P) &=&   m^{-1}(Ns p)^m(p+q+Ns)^{-m}\\
  \nonumber
  S_N(\bx,P) &=&\{m(m-1)\}^{-1}\left[(m-1)(t_N+1)^{a-m}\{(t_N+Ns)(t_N+p)\}^{m}(t_N+p+q+Ns)^{-m} \right.\\
  \label{eq:suffnegbinomform}
  &&{} - \left. mt_N^{a-m+1}\{(t_N+Ns-1)(t_N+p-1)\}^{m-1}(t_N+p+q+Ns-1)^{-m+1}\right].
\end{eqnarray}
 
\subsubsection{Prequential}
\label{sec:preqnegbinom}

Now suppose we have already observed $\bX^{n-1} = \bx^{n-1}$.  The
posterior distribution of $\Theta$ is
\begin{displaymath}
  \Theta \cd \bX^{n-1} = \bx^{n-1} \sim \beta\left\{p+t_{n-1},q+(n-1)s\right\}. 
\end{displaymath}

So the predictive density of $X_n$, given the previous
observations $\bX^{n-1} = \bx^{n-1}$, is obtained from
\eqref{eq:negbinomform0} and \eqref{eq:negbinomform} on replacing $x$ with $x_n$,
$p$ with $p+t_{n-1}$, and $q$ with $q+(n-1)s$.  The
incremental contribution to the prequential score is thus given by:

\begin{eqnarray}
  \label{eq:preqnegbinomform0}
  S^*_n(0,P) &=&  m^{-1}s^m (p+t_{n-1})^m(p+q+t_{n-1}+ns)^{-m}\\
  \nonumber
  S^*_n(x_n,P) &=&\{m(m-1)\}^{-1}\left[(m-1)(x_n+1)^{a-m}\left\{(x_n+s)(p+t_n)\right\}^m(p+q+t_n+ns)^{-m} \right.\\
  \label{eq:preqnegbinomform}
  &&{} - \left. mx_n^{a-m+1}\left\{(x_n+s-1)(p+t_n-1)\right\}^{m-1}(p+q+t_n+ns-1)^{-m+1}\right].
\end{eqnarray}
The total prequential score is obtained by summing this from $n=1$ to
$N$.

\subsection{Improper prior}

The usual improper prior is the formal limit with $p, q
\downarrow 0$.  In this case \eqref{eq:suffnegbinomform0} and \eqref{eq:suffnegbinomform} become:

\begin{eqnarray}
  \label{eq:impsuffnegbinomform0}
  S_N(\bo,P) &=&   0\\
  \nonumber
  S(x,P) &=&\{m(m-1)\}^{-1}\left\{(m-1)(t_N+1)^{a-m}t_N^{m}\right.\\
  \label{eq:impsuffnegbinomform}
  &&{} - \left. mt_N^{a-m+1}(t_N-1)^{m-1}\right\}.
\end{eqnarray}

The score is well-defined even when all observations are $0$, in which
case the posterior is improper.

For the prequential version, we obtain, from \eqref{eq:preqnegbinomform0} and
\eqref{eq:preqnegbinomform}:

\begin{eqnarray}
  \label{eq:preqnegbinomform01}
  S^*_n(0,P) &=&  m^{-1}s^m t_{n-1}^m(t_{n-1}+ns)^{-m}\\
  \nonumber
  S^*_n(x_n,P) &=&\{m(m-1)\}^{-1}\left[(m-1)(x_n+1)^{a-m}(x_n+s)^mt_n^m(t_n+ns)^{-m} \right.\\
  \nonumber
             &&{} - \left. mx_n^{a-m+1}(x_n+s-1)^{m-1}(t_n-1)^{m-1}(t_n +ns -1)^{-m+1}\right].\\
  \label{eq:preqnegbinomform1}\quad
\end{eqnarray}
The total prequential score is obtained by summing this from $n=1$ to
$N$.

Again, similar expressions can be found using the improper Jeffreys
prior, which has $p \downarrow 0$, $q = 1/2$.

\section{Simulations}
\label{sec:sim}

We generated observations from either the Poisson distribution
\eqref{eq:pois} with $k=1$, $\lambda = 10$, or the Negative Binomial
distribution \eqref{eq:negbin} with $s=81$, $\theta = 0.1$.  These
both have variance $10$, the former having mean $10$, and the latter
mean $9$.  We used, as the scoring rule, the special case of
\eqref{eq:sform} having $a=m=2$, namely
\begin{displaymath}
  S(x,P) =
  \frac 1 2 (x+1)^2\left\{\frac{p(x+1)}{p(x)}\right\}^2 - 
  x^2\left\{\frac{p(x)}{p(x-1)}\right\}\delta(x > 0).
\end{displaymath}

\begin{figure}[p]
  \centering
    \includegraphics[width=\linewidth]
    {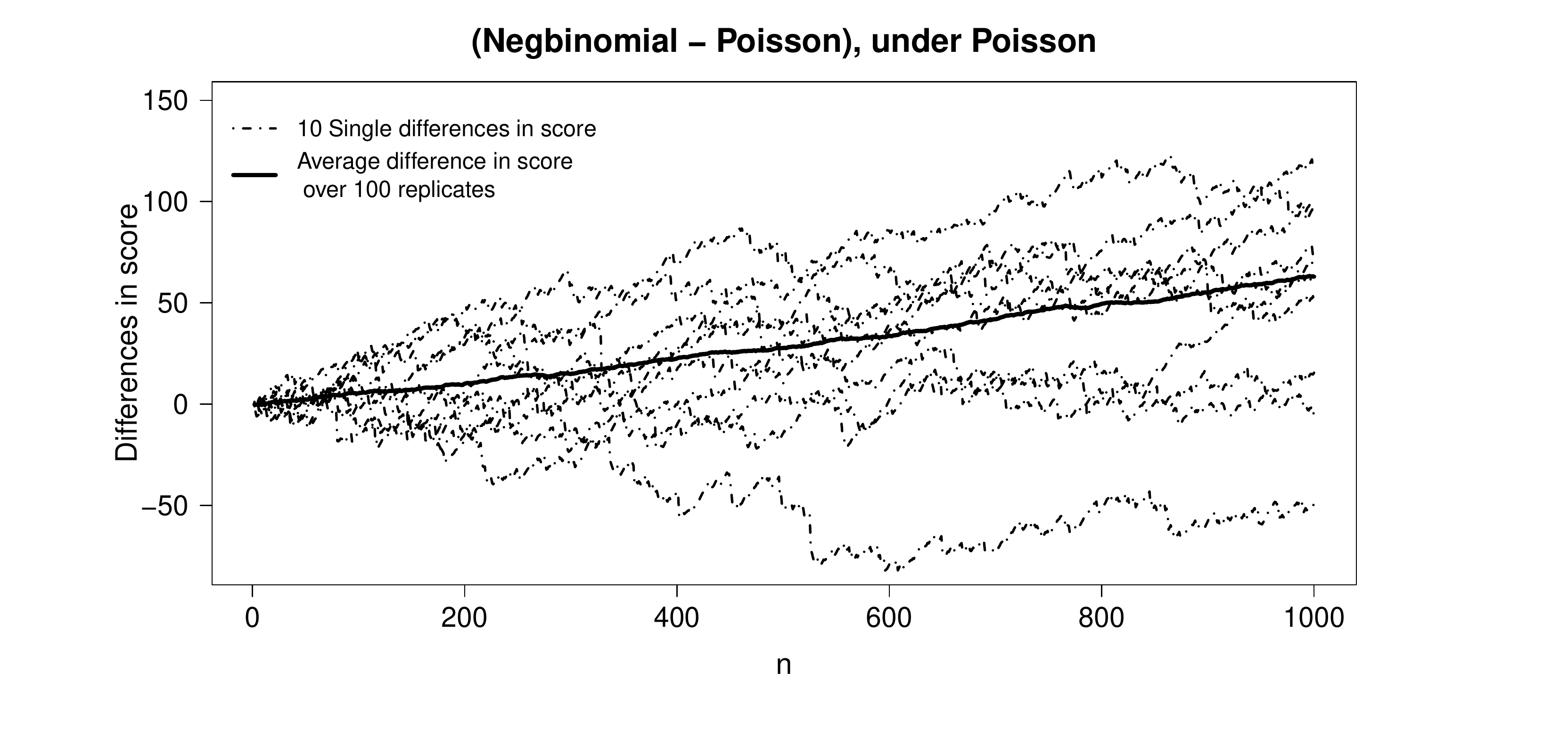}
    \caption{Data from Poisson distribution ${\cal P}(10)$}
    \label{fig:poissfig}
  \end{figure}

\begin{figure}[p]
  \centering
    \includegraphics[width=\linewidth]
    {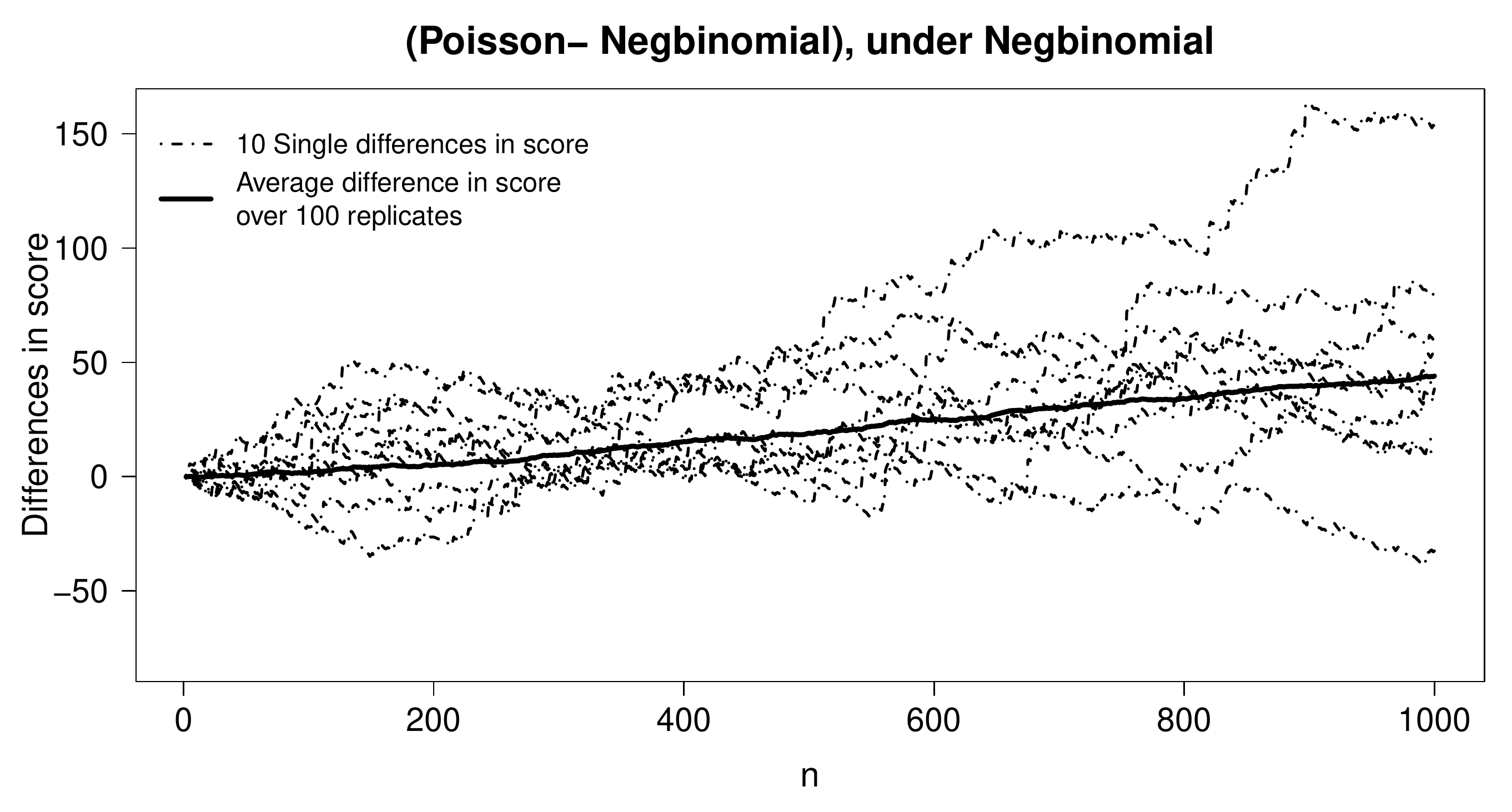}
    \caption{Data from Negative Binomial distribution ${\cal NB}(81;0.1)$}
    \label{fig:nbfig}
  \end{figure}

For each generating distribution we computed the excess of the
cumulative prequential score for the wrong model over that for the
correct model.  These differences are shown, as a function of
increasing data, in Figures~\ref{fig:poissfig} and \ref{fig:nbfig}
respectively.  Each figure displays 10 sample sequences generated from
the indicated distribution, as well as the average taken over a sample
Areof 100 sequences.

In each case we see a clear linear upward trend, supporting the
expectation of consistent model selection, although even with 1000
observations there is a non-negligible probability of a negative
value, giving a misleading preference for the wrong model.

\section{Conclusions}
\label{sec:conc}

We have extended the Bayesian model selection methodology of
\citet{apd/mm:modsel} to apply to problems with discrete data.  We have
conducted a simulation study to compare Poisson and Negative Binomial
distributions.  The results suggest that the method will consistently
select the correct model as the number of data points increases.


\section*{Acknowledgements}
\noindent Philip Dawid's research was supported through an Emeritus
Fellowship from the Leverhulme Trust.

\end{document}